\documentstyle[11pt]{article}
\newcommand{ \newsection}[1]{ \setcounter{equation}{0} \section{ #1} }

\def\keywords{ \if@twocolumn
\section*{Keywords}
\else \small
\begin{center}
{ \bf Keywords\vspace{-.5em}\vspace{0pt}}
\end{center}
\center
\fi}
\def\endkeywords{ \if@twocolumn\else\endcenter\fi}
\begin{document}
\begin{center}
{\large  \bf QUINTIC SPLINE SOLUTIONS OF FOURTH ORDER } \\
{ \large \bf BOUNDARY-VALUE PROBLEMS} \\
Shahid S. Siddiqi  \\
And \\
 Ghazala Akram \\
 Department of Mathematics, Punjab University,\\
Lahore, Postcode 54590 Pakistan. \\
\end{center}
\ \ \ \\
\begin{center}
\begin{minipage}{5.0in}
{\small Dedicated to the memory of Dr. M. Rafique} \\
\ \ \ \\
\\
{\bf KEYWORDS:}
{\small
Quintic spline; Boundary Value Problem; End conditions.} \\
  \\
 AMS Classification: 65L10
 \begin{abstract}
\noindent In this paper Quintic Spline is defined for the
numerical solutions of the fourth order linear special case
Boundary Value Problems. End conditions are also derived to
complete the definition of spline.The algorithm developed
approximates the solutions, and their higher order derivatives of
differential equations. Numerical illustrations are tabulated to
demonstrate the practical usefulness of method.
\end{abstract}
\end{minipage}
\end{center}
\ \ \\
\newsection{Introduction}
  Spline functions are used in many areas such as interpolation, data fitting,
   numerical solution of ordinary and partial differential equations. Spline functions
   are also used in curve and surface designing. \\
    Riaz A. Usmani [1], considered the fourth order boundary value problem
  to be the problem of bending a rectangular clamped beam of length l resting on an elastic
 foundation. The vertical deflection w of the beam satisfies the
 system
\begin{equation}
  \left[ L \ + \  \left( \frac{K}{D} \right ) \right ] w = D^{-1}q(x)  \  \ where \ \  L\equiv \frac{d^{4}}{dx^{4}}
\end{equation}
\begin{equation}
  w(0) = w(L) = w'(0) = w'(L) = 0
\end{equation}
where D is the flexural rigidity of the beam, and k is the spring
constant of the elastic foundation and the load q(x) acts
vertically downwards per unit length of the beam. The detail of
the mechanical interpretation of (1.1) belongs to a general class
of boundary value problems of the form \\
\newpage
 \pagestyle{myheadings}
\begin{equation}
\left. \begin{array}{c}
  \left( \frac{d^{4}}{dx^{4}}\ + \ f(x) \right ) y(x) \  = \   g(x), \   x \ \in \
  [a,b] \\
   y(a) = \alpha_{0},\  y(b) \  = \ \alpha_{1} \\
  y'(a) = \beta_{0},\ y'(b) \ = \  \beta_{1}
  \end{array} \right\}
\end{equation}
where $ \alpha_{i}, \ \beta_{i} \ ;$ i= 0,1 are finite real
constants and the functions f(x) and g(x) are continuous on [a,b].
The analytic solution of (1.3) for special choices of f(x) and
g(x) are easily obtained, but for arbitrary choices, the analytic
solution cannot be determined.\\
 Numerical methods for obtaining an approximation to y(x) are introduced. Usmani [1] derived
numerical techniques of order 2, 4 and 6 for solution of a fourth
order linear boundary value problem. Usmani [2] derived cubic,
quartic, quintic and sextic spline solution of nonlinear boundary
value problems. Usmani and Manabu Sakai [3] developed a quartic
spline for the approximation of the solution of certain two point
boundary value problems involving third order linear differential
equation.\\
 N.Papamichael and A.J.Worsey [4] derived end
conditions for cubic spline interpolation at equally spaced knots.
N.Papamichael and A.J.Worsey [5] have developed a cubic spline
method, similar to that proposed by Daniel and Swartz [6] for
second order problems. In this paper a quintic spline method is
described for the solution of (1.3). The end conditions for
quintic spline interpolation, at equally spaced knots are derived
which uniformly converge on [a,b] to O($h^6$), which is discussed
in the next section.
%
%
\newsection{Quintic Spline}

Let Q be a quintic spline defined on $[a,b]$ with equally spaced
knots
\begin{equation}
  x_{i} = a \ + \ ih ; \ \  i=0,1,2,...,k
\end{equation}
where
\begin{equation} h = \frac{b-a}{k}
\end{equation}
Moreover for i=0,1,2,\ldots, \ k, taking
\begin{equation}
       Q(x_{i}) = y_{i} \  , \  Q^{(1)}(x_{i}) = m_{i} ,
\end{equation}
\begin{equation}
    Q^{(2)}(x_{i}) = M_{i}  \ , \  Q^{(3)}(x_{i}) = n_{i},
\end{equation}
and
\begin{equation}
      Q^{(4)}(x_{i}) = N_{i}
\end{equation}
Also, let y(x) be the exact solution of the system (1.3) and
$y_{i}$ be an approximation to $ y(x_{i})$, obtained by the
quintic spline $Q(x_{i}).$ It may be noted that the $ Q_{i}(x), \
i=1,2,3,\ldots,k $ can be defined on the interval
$[x_{i-1},x_{i}]$, integrating
\begin{equation}
  {Q^{4}}_{i}(x) = \frac{1}{h}[N_{i-1}(x_{i} - x) + N_{i}(x-x_{i-1})]
\end{equation}
four times w.r.t. x, which gives
\begin{equation}
 Q_{i}(x) = \frac{1}{120 h}[N_{i-1}(x_{i} - x)^{5} \ + \
         N_{i}(x-x_{i-1})^{5}] \ + \ \frac{A x^{3}}{6} \ + \ \frac{B
         x^{2}}{2} \ + \ C x  \ + \ D
 \end{equation}
 To calculate the constants of integrations, the following
conditions are used.
\begin{eqnarray}
 Q_{i}(x_{i}) = y_{i} \ ; & \  {Q^{2}}_{i}(x_{i}) = M_{i}  \nonumber \\
  Q_{i}(x_{i-1}) = y_{i-1} \  ; & \  {Q^{2}}_{i}(x_{i-1}) = M_{i-1}
\end{eqnarray}
The identities of quintic splines for the solution of (1.3) can be
written as
\begin{eqnarray}
    m_{i-2}\ + \ 26m_{i-1} \ + \ 66m_{i} \ + \ 26 m_{i+1} \ + \
  m_{i+2}  &
  =&
  \frac{5}{h} \left[ -y_{i-2}-10 y_{i-1}+10 y_{i+1}+y_{i+2} \right ] \ , \nonumber \\
  & &  \ i=2,3,...,k-2
\end{eqnarray}
\begin{eqnarray}
  M_{i-2}\ + \ 26M_{i-1} \ + \ 66M_{i} \ + \ 26 M_{i+1} \ + \ M_{i+2}
   & = &
  \frac{20}{h^{2}}[ y_{i-2}+2 y_{i-1}-6y_{i}+2 y_{i+1}+y_{i+2}] \ ,
  \nonumber \\
  &  & i=2,3,...,k-2
\end{eqnarray}
\begin{eqnarray}
  n_{i-2}\ + \ 26n_{i-1} \ + \ 66n_{i} \ + \ 26 n_{i+1} \ + \ n_{i+2} &
  =&
  \frac{60}{h^{3}}[-y_{i-2}+2y_{i-1}-2y_{i+1}+y_{i+2}] \ , \nonumber
  \\
   & & i=2,3,...,k-2
\end{eqnarray}
\begin{eqnarray}
    N_{i-2}\ + \ 26N_{i-1} \ + \ 66N_{i} \ + \ 26 N_{i+1} \ + \
   N_{i+2}   & = &
  \frac{120}{h^{4}}[ y_{i-2}-4y_{i-1}+6y_{i}-4y_{i+1}+y_{i+2}] \ ,  \nonumber \\
   & &  i=2,3,...,k-2
  \end{eqnarray}
\begin{eqnarray}
  N_{i-1}\ + \ 4N_{i} \ + \ N_{i+1}- \frac{6}{h^{2}}[M_{i-1}-2M_{i}+M_{i+1}]& = &
  0 \ ,
  \nonumber \\
     i &= & 1,2,3,...,k-1
\end{eqnarray}
\begin{eqnarray}
 60h
  \{
     m_{i-1}\ + \ 2m_{i} \ + \ m_{i+1}
 \}
  - h^{3}
  \{
  3n_{i-1} \ +
 \ 14n_{i}\ + \ 3n_{i+1}
 \} &
 = &
  120\{- y_{i-1}+ y_{i+1}\} \ ,
  \nonumber \\
    i=1,2,3,...,k-1  & &
\end{eqnarray}
\begin{eqnarray}
  8h\{m_{i+1}\ - \  m_{i-1}\}-h^{2}\{M_{i-1} \ - \  6M_{i} \ + \ M_{i+1}\}&
  =&
  20\{y_{i-1}-2y_{i}+y_{i+1}\} \ ,
  \nonumber \\
     i & = & 1,2,3,...,k-1
\end{eqnarray}
\begin{eqnarray}
  m_{i} & = & \frac{h}{6} \{2 M_{i}\ + \  M_{i-1}\} \ -  \frac{h^{3}}{360}\{8N_{i} \ + \  7N_{i-1}\}
  \ + \ \frac{1}{h}\{y_{i}-y_{i-1}\} \ ,
  \nonumber \\
    &  & i=1,2,3,...,k
\end{eqnarray}
\begin{eqnarray}
  m_{i}& = & - \frac{h}{6} \{2 M_{i}\ + \ M_{i+1}\}+  \frac{h^{3}}{360} { 8N_{i} \ + \  7N_{i+1}}
  \ + \ \frac{1}{h}\{y_{i+1}-y_{i}\} \ ,
  \nonumber \\
     i & = & 0,1,,...,k-1
\end{eqnarray}
\begin{eqnarray}
  m_{i} & = & - \frac{h^{2}}{120} \{n_{i-1} \ + \  18n_{i} \ + \ n_{i+1}\}+
  \frac{1}{2h}\{y_{i+1}-y_{i-1}\} \ ,
  \nonumber \\
    i & = &1,2,3,...,k-1
\end{eqnarray}
\begin{eqnarray}
  M_{i} & =& - \frac{h^{2}}{120} \{N_{i-1} \ + \ 8N_{i} \  + \ N_{i+1}\}+ \frac{1}{h^{2}}
  \{y_{i-1}-2y_{i}+ y_{i+1}\} \ ,
  \nonumber \\
     i & = & 1,2,3,...,k-1
\end{eqnarray}
\begin{eqnarray}
  M_{i} & = &  \frac{1}{32h}  \{m_{i-2} \  + \ 32m_{i-1} \ - \ 32m_{i+1} \ -  \
  m_{i+2}\}  \nonumber \\
            &             & +  \frac{5}{32h^{2}}\{y_{i-2}+
            16y_{i-1}-34y_{i}+16y_{i+1}+y_{i+2}\} \ ,
  \nonumber \\
   &  & i=2,3,...,k-2
\end{eqnarray}
\begin{eqnarray}
  N_{i} = - \frac{3}{2h^{2}}  \{M_{i-1}+&18M_{i} \ + \ M_{i+1}\}
        +  \frac{30}{h^{4}}\{y_{i-1}-2y_{i}+y_{i+1}\} \ ,
  \nonumber \\
    & i=1,2,3,...,k-1
\end{eqnarray}
The relations (2.9)--(2.21) can be derived from the results of
Albasiny and Hoskins(1971), Fyfe(1971), Ahlberg, Nilson and
Walsh(1966) and Sakai(1970) discussed by N. Papamichael [7].
 The uniqueness of Q can be established showing that any of the
  four $(k+1) \times (k+1)$ linear systems, obtained, using one of the relations
   (2.9), (2.10), (2.11) and (2.12)together with the four end conditions.
   In this paper the linear
system corresponding to (2.12) has been chosen and has a unique
solution for $ N_{i}s \ , \ i=0,1,...,k$. Equation (2.13) and
(2.19), or (2.19) and (2.21) give the parameters $M_{i}s \ , \
i=0,1,...,k$. Consider the system (2.12)
\begin{eqnarray}
       N_{i-2}\ + \ 26N_{i-1} \ + \ 66N_{i}   \  +   \ 26 N_{i+1} \ + \
  N_{i+2} &
    = &
  \frac{120}{h^{4}} \ [y_{i-2}-4y_{i-1}+6y_{i}-4y_{i+1}+y_{i+2}] \ ,
  \nonumber \\
 & &      i=2,3,...,k-2
\end{eqnarray}
where
\begin{equation}
N_{i} = -f_{i}y_{i} + g_{i}
\end{equation}
The above system gives $(k-3) $ linear algebraic equations in the
$(k-1)$ unknowns $(y_{i}, \ i=1,2,...,k-1) $. Two more equations
are needed to have complete solution of $y_{i}s$ which are derived
in section 3.\\
 \textbf{It may be noted that Papamichael [5]
needed two consistency systems for the solution of boundary value
problem (1.3) but the method developed in this paper, needs only
one such system.}\\
 Taking  forward difference operator $ E \ = \ e^{hD}$, equation
 (2.12) can be rewritten as
 \begin{eqnarray}
 (E^{-2}  +  26 E^{-1}  +  66I  +  26 E^{1}  +
   E^{2} )N_{i}  & = &
  \frac{120}{h^{4}} \ [ \ E^{-2}  -  4E^{-1}  +  6I-  4E^{1}  +  E^{2}] \ y_{i} \ ,  \nonumber \\
   & &  i=2,3,...,k-2
  \end{eqnarray}
\begin{eqnarray}
N_{i}& = & \frac{120}{h^{4}} \ ( \ e^{-2hD} \ + \ 26 \ e^{-hD} \ +
           \ 66 \ I  \ + \ 26 \ e^{hD} \ + \ e^{2hD} )^{-1} \nonumber \\
     &  & ( e^{-2hD} \ - \ 4 \ e^{-hD} \ + \ 6 \ I \ - \ 4 \ e^{hD} \ + \ e^{2hD}) \ y_{i} \ ,  \nonumber \\
    i & = & 2,3,...,k-2
   \end{eqnarray}
  Expanding the R.H.S. of above, in terms of power series and dividing give,
   \begin{equation}
   N_{i} \ = \ y_{i}^{(4)} \ - \ \frac{h^2}{12} \ y_{i}^{(6)} \ + \
   \frac{h^4}{240} \ y_{i}^{(8)} \ + \ O(h^{6})
   \end{equation}
 Papamichael [7] proved the following lemma to determine the end conditions
 in terms of first derivative of quintic spline \\
\\
{\large {\bf Lemma 1 }}

 Let $ \lambda_{i} \ = \ m_{i} \ - \ y_{i}^{(1)}$. If
$ y \in C^{7}[a,b] $ then
\begin{equation}
  \lambda_{i-2} +26 \lambda_{i-1}+66 \lambda_{i}+26
  \lambda_{i+1}+\lambda_{i+2}= \beta_{i} \ , \  i=2,3,...,k-2
\end{equation}
where \\ $$ |\beta_{i}|  \ \leq \ \frac{11}{21} \  h^{6} \  \|y^{(7)}\| \ , \ \ i=2,3,...,k-2 $$ \\
Using the above lemma along with equation(2.12) and with Taylor
series expansion about the point $x_{i}$, the following lemma can
easily be proved, for the case discussed
 in this paper. \\
 \\
{\large {\bf Lemma 2 }}

 Let
\begin{equation}
\lambda_{i} = y^{(4)}_{i} \ - \ \frac{h^{2}}{12} \  y^{(6)}_{i} +
 \ \frac{h^{4}}{240} \  y^{(8)}_{i} -N_{i}
  \end{equation}
If $y \in C^{10}[a,b]$ then
\begin{equation}
  \lambda_{i-2} +26 \lambda_{i-1}+66 \lambda_{i}+26
  \lambda_{i+1}+\lambda_{i+2}= \beta_{i} \ , \  i=2,3,...,k-2
\end{equation}
where \\ $$ |\beta_{i}|  \ \leq \ \frac{2665920}{10!} \  h^{6} \  \|y^{(10)}\| \ , \ \ i=2,3,...,k-2 $$ \\
Finally the required end conditions are derived in the following
section.
\newsection{End Conditions}
Consider end conditions of the form
\begin{equation}
  N_{0}\ + \ \alpha N_{1}\ + \ \beta N_{2} \ + \  \gamma N_{3} \ + \ N_{4}
  \\ = \frac{1}{h^{4}} \left[ \  \sum_{i=0}^{3}a_{i}y_{i}\ + \
  bhy_{0}^{(1)}\ + \ h^{4} \ \sum_{i=0}^{4}c_{i}y_{i}^{(4)} \  \right ]
  \end{equation}
  and
\begin{equation}
   N_{k}\ + \ \alpha N_{k-1}\ + \ \beta N_{k-2} \ + \ \gamma N_{k-3} \ +
  \ N_{k-4} \\
  =   \frac{1}{h^{4}} \left[  \ \sum_{i=0}^{3}a_{i}y_{k-i} \ + \ bhy_{k}^{(1)}
    \ +  h^{4} \ \sum_{i=0}^{4}c_{i}y_{k-i}^{(4)} \  \right ]
\end{equation}
where the scalars $ \alpha, \ \beta, \ \gamma $ , $a_{i}$,
i=0,1,2,3, b and $c_{i}$, i=0,1,...,4 can be determined with the
assumption that Q exists uniquely and
 $$ \|Q^{r} \ - \ y^{r} \| \ = \ O(h^{6-r}) \ , \ \  r= 0,1,...,5 $$
 For this let $$\lambda_{i} \ = \ y^{(4)}_{i}-\frac{h^{2}}{12} \
y^{(6)}_{i} \ + \  \frac{h^{4}}{240} \  y^{(8)}_{i} \ - \ N_{i} \
,\ \ i=0,1,...,k.$$ Since it is supposed that $ y \in C^{10}[a,b]
$, therefore the equations (2.12), (3.1) and (3.2) give,
\begin{equation}
 \lambda_{0} \ + \ \alpha\lambda_{1} \ + \ \beta\lambda_{2}\ + \
 \gamma\lambda_{3}\ + \ \lambda_{4} \ = \ \beta_{1}
\end{equation}
\begin{eqnarray}
  \lambda_{i-2} \ + \ 26\lambda_{i-1} \ + \ 66\lambda_{i}\ + \
 26\lambda_{i+1}\ + \ \lambda_{i+2} \ = \ \beta_{i} \ ,  \nonumber \\  i=2,3,...,k-2
\end{eqnarray}
\begin{equation}
  \lambda_{k} \ + \ \alpha\lambda_{k-1} \ + \ \beta\lambda_{k-2}\ + \
 \gamma\lambda_{k-3}\ + \ \lambda_{k-4} \ = \ \beta_{k-1}
\end{equation}
where
\begin{eqnarray}
  \beta_{1}\ & = & \ \frac{1}{h^{4}} \  \left[ - \sum_{i=0}^{3}a_{i}y_{i}\ - \
  bhy_{0}^{(1)}\ - \ \sum_{i=0}^{4}c_{i}y_{i}^{(4)} \ + \ h^{4} \ ( \ y_{0}^{(4)} \ + \ \alpha y_{1}^{(4)}
    \ + \ \beta y_{2}^{(4)} \ + \ \gamma y_{3}^{(4)} \ + \ y_{4}^{(4)} \ )
    \right.
    \nonumber \\
    & &  \ - \ \frac{h^{6}}{12} \ ( \ y_{0}^{(6)} \ + \ \alpha y_{1}^{(6)}
    \ + \ \beta y_{2}^{(6)} \ + \ \gamma y_{3}^{(6)} \ + \ y_{4}^{(6)} \ )
    \nonumber \\
    & &  \left.  \ + \ \frac{h^{8}}{240} \ ( \ y_{0}^{(8)} \ + \ \alpha y_{1}^{(8)}
    \ + \ \beta y_{2}^{(8)} \ + \ \gamma y_{3}^{(8)} \ + \ y_{4}^{(8)} \ ) \right]
\end{eqnarray}
and
\begin{eqnarray}
 \beta_{k-1}\ & = & \ \frac{1}{h^{4}} \ \left[ - \sum_{i=0}^{3}a_{i}y_{k-i}\ - \
  bhy_{k}^{(1)}\ - \ \sum_{i=0}^{4}c_{i}y_{k-i}^{(4)} \ + \ h^{4} \ ( \ y_{k}^{(4)} \ + \ \alpha y_{k-1}^{(4)}
    \ + \ \beta y_{k-2}^{(4)} \right. \ \nonumber \\
    & & + \ \gamma y_{k-3}^{(4)} \
     \ + \ y_{k-4}^{(4)}) \ - \ \frac{h^{6}}{12} \ ( \ y_{k}^{(6)} \ + \ \alpha y_{k-1}^{(6)}
    \ + \ \beta y_{k-2}^{(6)} \ + \ \gamma y_{k-3}^{(6)} \ + \ y_{k-4}^{(6)} \ )
    \nonumber \\
    & &  \left. \ + \ \frac{h^{8}}{240} \ ( \ y_{k}^{(8)} \ + \ \alpha y_{k-1}^{(8)}
    \ + \ \beta y_{k-2}^{(8)} \ + \ \gamma y_{k-3}^{(8)} \ + \ y_{k-4}^{(8)} \ ) \right]
\end{eqnarray}
and, from lemma 2
\begin{equation}
  |\beta_{i}|  \ \leq \ \frac{2665920}{10!} \  h^{6} \  \|y^{(10)}\| \ , \ \ i=2,3,...,k-2
\end{equation}
Following Papamichael [7], the required end conditions may be
written as
\begin{eqnarray}
  N_{0}\ + \  N_{1}\ \ & = & \frac{1}{h^{4}} \ \left[  \ \frac{220}{3}y_{0}\
   - \ 120 y_{1} \ + \ 60 y_{2} \ - \frac{40}{3}y_{3} \ + \
   40hy_{0}^{(1)} \right.
   \nonumber \\
    & &  \ + \ h^{4} \left(  \ \frac{2519}{2520} \ y_{0}^{(4)}\ + \ \frac{11822}{1260} \ y_{1}^{(4)}
   \ + \ \frac{223}{210} \ y_{2}^{(4)} \ + \
   \frac{176}{180} \ y_{3}^{(4)}\right.
   \nonumber \\
    & &  \left. \left. \ + \ \frac{1769}{2520} \ y_{4}^{(4)} \right) \right ]
\end{eqnarray}
and
\begin{eqnarray}
  N_{k}\ + \  N_{k-1}\ \ & = & \frac{1}{h^{4}} \ \left[ \ \frac{220}{3} \ y_{k}\
   - \ 120 y_{k-1} \ + \ 60 y_{k-2} \ - \frac{40}{3} \ y_{k-3} \ + \ 40hy_{k}^{(1)}
    \right. \nonumber \\
  & & \ + \ h^{4} \left( \frac{2519}{2520} \ y_{k}^{(4)}\ + \ \frac{11822}{1260} \ y_{k-1}^{(4)}
   \ + \ \frac{223}{210} \ y_{k-2}^{(4)} \ + \ \frac{176}{180} \ y_{k-3}^{(4)}
    \right.  \nonumber \\
   & &
    \left. \left. \ + \ \frac{1769}{2520} \ y_{k-4}^{(4)} \right) \right]
\end{eqnarray}
The quintic spline solution of the system (1.3) is defined in the
next section.
\newsection{Quintic Spline Solution} The quintic spline solution of
(1.3) is based on the linear equations (2.12), (3.9) and (3.10).
Let \textbf{Y}$ \ = \ y_{i}$, \ \textbf{C}$ \ = \ c_{i}$,
\textbf{e} $ \ = \ e_{i}$. Then the parameters $y_{i}$ of Q
satisfy the following matrix equation
$$
(\textbf{A} \ + \ h^{4} \ \textbf{B}\textbf{F})\textbf{Y} \ = \
\textbf{C}\ + \ \textbf{e}
$$
where \textbf{Y}, \textbf{C}, \textbf{e} are $(k-1)$ dimensional
column vectors and \textbf{A}, \textbf{B}, \textbf{F} are $(k-1) \
\times  \ (k-1)$ matrices, where \\
\\
\\
 $$ A= \left[ \begin{array}{rrrrrrrr}
     9      & \frac{-9}{2} & 1       &      &        &       &      &   \\
     &  & & & & & & \\
    -4      & 6           & -4      & 1    &        &       &      &   \\
    &  & & & & & & \\
     1      & -4          & 6       & -4    & 1     &       &      &  \\
     &  & & & & & & \\
            &   \ddots         & \ddots       & \ddots     & \ddots      & \ddots      &      & \\
              &  & & & & & & \\
            &             &         &1      &-4     &6      &-4    & 1 \\
            &  & & & & & & \\
            &             &         &       &1      &-4     &6      &-4 \\
            &  & & & & & & \\
            &             &         &       &       &1      & \frac{-9}{2} & 9
      \end{array} \right] $$
      \\
    $$  B= \left[ \begin{array}{rrrrrrrr}
    \frac{ 31686}{5040}& \frac{669}{8400} & \frac{528}{7200}&\frac{5307}{100800}& & &  &   \\
    &  & & & & & & \\
    \frac{13}{60}      & \frac{11}{20}    & \frac{13}{60}   &\frac{1}{120}      & & &  &   \\
     &  & & & & & & \\
    \frac{1}{120}      & \frac{13}{60}    & \frac{11}{20}   & \frac{13}{60}    & \frac{1}{120}&&&\\
     &  & & & & & & \\
            & \ddots        & \ddots               & \ddots     & \ddots     &\ddots     &      & \\
     &  & & & & & & \\
            &          &         &\frac{1}{120}      & \frac{13}{60}    & \frac{11}{20}   & \frac{13}{60}    & \frac{1}{120} \\
     &  & & & & & & \\
            &          &         &       &\frac{1}{120}      & \frac{13}{60}    & \frac{11}{20}   &\frac{13}{60} \\
     &  & & & & & &  \\
            &         &         &        &\frac{5307}{100800}& \frac{528}{7200} & \frac{669}{8400}&\frac{31686}{5040}
      \end{array} \right]  $$
  \\
      $ \textbf{F} = diag(f_{i}) $
\begin{eqnarray}
c_{1} & = & \frac{11}{2} \ y_{0}\ + \ 3hy_{0}^{(1)} \ + \
            \frac{h^{4}}{280} \ \left( \frac{-3}{360} \ g_{0} \ + \
            \frac{3}{360} \ f_{0} \ y_{0} \ + \ \frac{31686}{180} \ g_{1}
            \right.
             \nonumber \\
      &   &  \left. \ + \ \frac{669}{30} \ g_{2} \ + \ \frac{3696}{180} \ g_{3} \ + \
\frac{5307}{360} \ g_{4} \right)
\end{eqnarray}
\begin{equation}
c_{2} = \frac{h^{4}}{120} \ ( \ g_{0}-f_{0} \ y_{0}+26 \
         g_{1}+66 \ g_{2}+ 26 \ g_{3}+g_{4} \ )-y_{0}
\end{equation}
\begin{eqnarray}
c_{i} & = & \frac{h^{4}}{120} \ ( \ g_{i-2} \ + \ 26 \ g_{i-1} \
            + \ 66 \ g_{i} \ + \ 26 \ g_{i+1} \ + \ g_{i+2} \ )
             \nonumber \\
      &   &   i=3,4,...,k-3
\end{eqnarray}
\begin{equation}
c_{k-2} \  = \  \frac{h^{4}}{120} \
             ( \ g_{k} \ - \ f_{k} \ y_{k} \ + \ 26 \ g_{k-1} \ + \ 66 \ g_{k-2} \
             + \ 26 \ g_{k-3} \ + \ g_{k-4} \ ) \ - \ y_{k}
\end{equation}
\begin{eqnarray}
c_{k-1} & = & \frac{11}{2} \ y_{k}\ - \ 3hy_{k}^{(1)} \ + \
            \frac{h^{4}}{280} \ \left( \frac{-3}{360} \ g_{k} \ + \
            \frac{3}{360} \ f_{k} \ y_{k} \ + \ \frac{31686}{180} \
            g_{k-1} \right.
             \nonumber \\
      &   & \left.  \ + \ \frac{669}{30} \ g_{k-2} \ + \ \frac{3696}{180} \ g_{k-3} \ + \
            \frac{5307}{360} \ g_{k-4} \ \right)
\end{eqnarray}
 and $ \textbf{e}_{i} = O(h^{6}) \ \ i=1,...,k-1 $ \\
 \textit{\textbf{Extending the method for the solution of sixth,
eighth and higher order boundary value problems, is in process.}}
\\
 To implement the method for the quintic spline solution of the
boundary value problem (1.3), three examples are discussed in the
following section.\\
\newsection{Numerical Examples}
 In this section numerical technique discussed in section \textbf{4} is illustrated, by the
 following three boundary value problems of the type
 (1.3).\\
  \ \ \\
{\large {\bf Example 1 }} \\
 \ \ \\
 Consider the following boundary value problem
\begin{equation}
\left.
\begin{array}{l}
  y^{iv} \ + \ 4 y \ = \ 1,  \ \ \ y(-1) \ = \ y(1) \ = \ 0,  \\
  \  \   \  \\
 \ y^{'}(-1)=
             \ - \ y^{'}(1) \  = \ \frac{\sinh(2) \ - \ \sin(2) }{4(\cosh(2) \ + \
             \cos(2))}
\end{array}
\right\}
\end{equation}
The analytic solution of the above problem is \\
\begin{eqnarray}
 y(x) & = & 0.25[1-2[\sin(1) \ \sinh(1) \ \sin(x) \ \sinh(x)
  \nonumber
 \\
      & &  +  \cos(1) \ \cosh(1) \ \cos(x) \ \cosh(x) ]/(\cos(2) \ + \
      \cosh(2))] \nonumber
      \end{eqnarray}
      The results are summarized in Table \textbf{1}.
       \begin{table}[htp]
\caption{ Maximum absolute errors for Problem 5.1 in
 $ y_{i}^{( \mu )}, \ \mu = 0,1,2,3,4 $ }
 \begin{center}
\begin{tabular}{|l|l|l|l|l|l|} \hline
               &                              &
               &                              &       &      \\
h              &  $ y_{i}$  &
                  $ y^{(1)}_{i}$  &
                  $y^{(2)}_{i} $ &
                  $y^{(3)}_{i} $  &
                  $y^{(4)}_{i}$   \\
               &               &                                   &
                                                   &
                                                   &          \\ \hline
               &                              &
               &                              &       &      \\
$\frac{1}{8}$  & $ 2.1 \times {10}^{-3} $
               & $ 5.5 \times {10}^{-3} $
               & $ 9.7 \times {10}^{-3} $
               & $ 1.06 \times {10}^{-2} $
               & $ 8.5 \times {10}^{-3} $                       \\ \hline
               &                              &
               &                              &       &      \\
$\frac{1}{16}$   & $ 9.10 \times {10}^{-5} $
               & $ 3.34 \times {10}^{-4} $
               & $ 8.86 \times {10}^{-4} $
               & $ 1.4 \times {10}^{-3} $
               & $ 3.64 \times {10}^{-4} $                          \\ \hline
               &                              &
               &                              &       &      \\
$\frac{1}{32}$  & $ 5.12 \times {10}^{-6} $
               & $ 1.77 \times {10}^{-5} $
               & $ 1.20 \times {10}^{-4} $
               & $ 2.60 \times {10}^{-4} $
               & $2.05 \times {10}^{-5} $                             \\ \hline
               &                              &
               &                              &       &      \\
$\frac{1}{64}$ & $ 2.85 \times {10}^{-6} $
               & $ 4.50 \times {10}^{-6} $
               & $ 2.44 \times {10}^{-5} $
               & $ 6.03 \times {10}^{-5} $
               & $ 1.14 \times {10}^{-5} $                            \\ \hline
               &                              &
               &                              &       &      \\
 $\frac{1}{128}$ & $ 7.87 \times {10}^{-7} $
               & $ 1.18 \times {10}^{-6} $
               & $ 5.83 \times {10}^{-6} $
               & $ 1.47 \times {10}^{-5} $
               & $ 3.23 \times {10}^{-6} $        \\ \hline
               &                              &
               &                              &       &      \\
$\frac{1}{256} $& $ 1.98 \times {10}^{-7} $
               & $ 2.79 \times {10}^{-7} $
               & $ 1.44 \times {10}^{-6} $
               & $ 3.59 \times {10}^{-6} $
               & $ 8.73 \times {10}^{-7} $  \\ \hline
               &                              &
               &                              &       &      \\
 $\frac{1}{512}$ & $ 4.15 \times {10}^{-8} $
               & $ 4.04 \times {10}^{-8} $
               & $ 3.05 \times {10}^{-7} $
               & $ 6.81 \times {10}^{-7} $
               & $ 2.46 \times {10}^{-7} $    \\ \hline
               &                              &
               &                              &       &      \\
 $\frac{1}{1024}$ & $ 1.07 \times {10}^{-7} $
               & $ 1.91 \times {10}^{-7} $
               & $ 8.01 \times {10}^{-7} $
               & $ 2.12 \times {10}^{-6} $
               & $ 3.49 \times {10}^{-7} $    \\ \hline
\end{tabular}
\end{center}
\end{table}
\newpage
{\large {\bf Example 2}} \\

  Consider the following boundary value problem
\begin{equation}
\left.
\begin{array}{ll}
   & y^{iv} \ + \ x y \ = \ -(8 \ + \ 7x \ + \ x^{3}) e^{x},    \\
 y(0) \ = \ y(1) \ = \ 0,
 \ & y^{'}(0) \ = \ 1, \
             \ y^{'}(1) \ = \ -e
\end{array}
\right\}
\end{equation}
The analytic solution of the above differential system is
$$ y(x) \ = \ x(1-x) \  e^{x} $$
The observed maximum errors (in absolute value) associated with $
y_{i}^{( \mu )}, \ \mu = 0,1,2,3,4 $ , for the system 5.2 are
briefly summarized in Table 2.
\begin{table}[htp]
\caption{ Maximum absolute errors for Problem 5.2 in
 $ y_{i}^{( \mu )}, \ \mu = 0,1,2,3,4 $ }
 \begin{center}
\begin{tabular}{|l|l|l|l|l|l|} \hline
               &                              &
               &                              &       &      \\
h              &  $ y_{i}$  &
                  $ y^{(1)}_{i}$  &
                  $y^{(2)}_{i} $ &
                  $y^{(3)}_{i} $  &
                  $y^{(4)}_{i}$   \\
               &               &                                   &
                                                   &
                                                   &          \\ \hline
               &                              &
               &                              &       &      \\
$\frac{1}{8}$  & $ 5.3828 \times {10}^{-4} $
               & $ 1.9 \times {10}^{-3} $
               & $ 6.0 \times {10}^{-3} $
               & $ 3.8 \times {10}^{-2} $
               & $ 3.2768 \times {10}^{-4} $                       \\ \hline
               &                              &
               &                              &       &      \\
$\frac{1}{16}$   & $ 8.8114 \times {10}^{-5} $
               & $ 2.9045 \times {10}^{-4} $
               & $ 1.4 \times {10}^{-3} $
               & $ 9.9 \times {10}^{-3} $
               & $ 5.3779 \times {10}^{-5} $                          \\ \hline
               &                              &
               &                              &       &      \\
$\frac{1}{32}$  & $ 1.636 \times {10}^{-5} $
               & $ 5.1554 \times {10}^{-5} $
               & $ 3.3586 \times {10}^{-4} $
               & $ 2.6 \times {10}^{-3} $
               & $9.824\times {10}^{-6} $                             \\ \hline
               &                              &
               &                              &       &      \\
$\frac{1}{64}$ & $ 3.3894 \times {10}^{-6} $
               & $ 1.058 \times {10}^{-5} $
               & $ 8.2872 \times {10}^{-5} $
               & $ 6.6254 \times {10}^{-4} $
               & $ 1.991 \times {10}^{-6} $                            \\ \hline
               &                              &
               &                              &       &      \\
 $\frac{1}{128}$ & $ 7.5919 \times {10}^{-7} $
               & $ 2.388 \times {10}^{-6} $
               & $ 2.269 \times {10}^{-5} $
               & $ 1.6824 \times {10}^{-4} $
               & $ 4.3958 \times {10}^{-7} $        \\ \hline
               &                              &
               &                              &       &      \\
$\frac{1}{256} $& $ 1.7884 \times {10}^{-7} $
               & $ 5.671 \times {10}^{-7} $
               & $ 5.9302 \times {10}^{-6} $
               & $ 4.2401 \times {10}^{-5} $
               & $ 1.0264 \times {10}^{-7} $  \\ \hline
               &                              &
               &                              &       &      \\
 $\frac{1}{512}$ & $ 3.7364 \times {10}^{-8} $
               & $ 1.1843 \times {10}^{-7} $
               & $ 1.3068 \times {10}^{-6} $
               & $ 9.4838 \times {10}^{-6} $
               & $ 2.1304 \times {10}^{-8} $    \\ \hline
               &                              &
               &                              &       &      \\
 $\frac{1}{1024}$ & $ 6.3664 \times {10}^{-8} $
               & $ 2.001 \times {10}^{-7} $
               & $ 1.9244 \times {10}^{-6} $
               & $ 9.969 \times {10}^{-6} $
               & $ 3.5439 \times {10}^{-8} $    \\ \hline
\end{tabular}
\end{center}
\end{table}
\newpage
{\large {\bf Example 3}} \\

 Consider the differential system
\begin{equation}
\left.
\begin{array}{ll}
 & y^{iv} \ - y \ = \ - \ 4( 2 x \ cos(x) \ + 3 \ sin(x)) ,   \\
 y(0) \ = \ y(1) \ = \ 0 \ , \ &
 \ y^{'}(-1) \ = \ 2  \ sin(1) \ , \
             \ y^{'}(1) \ = \ 2 \  sin(1)
\end{array}
\right\}
\end{equation}
The analytic solution of the above system is
$$ y(x) \ = \ (x^{2}-1) \ sin(x) $$
The observed maximum errors (in absolute value) associated with $
y_{i}^{( \mu )}, \ \mu = 0,1,2,3,4 $ , for the system 5.3 are
briefly summarized in Table 3.
\begin{table}[htp]
\caption{ Maximum absolute errors for Problem 5.3 in
 $y_{i}^{( \mu )}, \ \mu = 0,1,2,3,4 $ }
 \begin{center}
\begin{tabular}{|l|l|l|l|l|l|} \hline
               &                              &
               &                              &      &       \\
  h            &  $ y_{i} $  &
                  $ y^{(1)}_{i} $ &
                  $ y^{(2)}_{i} $ &
                  $y^{(3)}_{i} $ &
                  $ y^{(4)}_{i} $  \\
               &               &                                   &
                                                   &
                                                  &          \\ \hline
                &                              &
               &                              &       &      \\
 $\frac{1}{8} $    & $ 1.4 \times {10}^{-3} $
               & $ 9.9 \times {10}^{-3} $
               & $ 4.3 \times {10}^{-2} $
               & $ 1.0 \times {10}^{-1} $
               & $ 1.4 \times {10}^{-3} $                       \\ \hline
               &                              &
               &                              &       &      \\
 $\frac{1}{16} $  & $ 1.42 \times {10}^{-4} $
               & $ 1.2 \times {10}^{-3} $
               & $ 7.0 \times {10}^{-3} $
               & $ 2.26 \times {10}^{-2} $
               & $ 1.42 \times {10}^{-4} $                          \\ \hline
               &                              &
               &                              &       &      \\
 $\frac{1}{32}$  & $ 8.18 \times {10}^{-6} $
               & $ 1.34 \times {10}^{-7} $
               & $ 1.2 \times {10}^{-3} $
               & $ 5.1 \times {10}^{-3} $
               & $ 8.81 \times {10}^{-6} $                             \\ \hline
               &                              &
               &                              &       &      \\
 $ \frac{1}{64} $ & $ 5.05 \times {10}^{-6} $
               & $ 1.95 \times {10}^{-5} $
               & $ 2.51 \times {10}^{-4} $
               & $ 1.2 \times {10}^{-3} $
               & $ 5.05 \times {10}^{-6} $                            \\ \hline
               &                              &
               &                              &       &      \\
 $ \frac{1}{128} $ & $ 1.6 \times {10}^{-6} $
               & $ 5.88 \times {10}^{-6} $
               & $ 5.68 \times {10}^{-5} $
               & $ 3.01 \times {10}^{-4} $
               & $ 1.61 \times {10}^{-6} $        \\ \hline
               &                              &
               &                              &       &      \\
 $\frac{1}{256} $ & $ 4.42 \times {10}^{-7} $
               & $ 1.58 \times {10}^{-6} $
               & $ 1.35 \times {10}^{-5} $
               & $ 7.48 \times {10}^{-5} $
               & $ 4.42 \times {10}^{-7} $  \\ \hline
               &                              &
               &                              &       &      \\
 $\frac{1}{512}$ & $ 1.15 \times {10}^{-7} $
               & $ 4.09 \times {10}^{-7} $
               & $ 3.31 \times {10}^{-6} $
               & $ 1.86 \times {10}^{-5} $
               & $ 1.15 \times {10}^{-7} $    \\ \hline
               &                              &
               &                              &       &      \\
 $ \frac{1}{1024}$ & $ 3.19 \times {10}^{-8} $
               & $ 1.05 \times {10}^{-7} $
               & $ 8.51 \times {10}^{-7} $
               & $ 5.01 \times {10}^{-6} $
               & $ 3.19 \times {10}^{-8} $    \\ \hline
\end{tabular}
\end{center}
\end{table}

\end{document}